\documentclass{amsart}
\usepackage[T1]{fontenc}
\usepackage{amsmath,amsthm,amssymb,amsfonts}
\newtheorem{theorem}{Theorem}

\newtheorem{corollary}[theorem]{Corollary}

\newtheorem{proposition}[theorem]{Proposition}
\newtheorem{remark}[theorem]{Remark}
\usepackage{enumitem}
\usepackage{amsthm}
\newtheorem*{theorem*}{Theorem A}
\newtheorem*{teorema}{Theorem B}
\newcommand{\Rp}{\mathbb{R}_+}
\newcommand{\al}{\alpha}

\begin{document}
\title[]{The weighted Bergman space on a sector and a degenerate parabolic equation}
\author{Marcos L\'{o}pez-Garc\'{\i}a}

\thanks{
The author was partially supported by project PAPIIT IN100919 of DGAPA-UNAM, and project A1-S-17475 of Conacyt, M\'exico.}
\email{marcos.lopez@im.unam.mx}
\address{
Instituto de Matem\'{a}ticas-Unidad Cuernavaca \\
   Universidad Nacional Aut\'{o}noma de M\'{e}xico\\
   Apdo. Postal 273-3, Cuernavaca Mor. CP 62251, M\'exico}

\maketitle
\begin{abstract}
In this work we solve a degenerate parabolic equation for the half line with Dirichlet boundary data, and use some results from the theory of Reproducing Kernel Hilbert Spaces to show that the null reachable space of this degenerate parabolic equation is a RKHS of analytic functions on a sector, whose reproducing kernel can be written in terms of the weighted Bergman kernel on the half plane $\mathbb{C}_+$. 
\end{abstract}

\section{Introduction}
Let $T>0$ fixed.  Consider the heat equation for the unit interval with Dirichlet boundary conditions,
\begin{eqnarray*}
\partial_t w-\partial_{xx} w=0, && 0<x<1, \, 0<t<T,  \notag  \\ 
w(0,t)=u_{\ell}(t), \quad w(1,t)=u_r(t), && 0<t<T ,\\
w(x,0)=0, &&  0<x<1.\notag
\end{eqnarray*}
In Control Theory of PDEs is an important issue to describe the so-called null reachable space, at time $T>0$, defined as follows
$$\mathcal{R}_T:=\{ w(\cdot,T): w \text{ is solution of the heat system with data }u_{\ell}, u_{r}\in L^2_{\mathbb{C}}(0,T)\}.$$

It is known that $\mathcal{R}_T$ does not depend on the time $T$, see \cite[Proposition 3.1]{kellay}. The problem is to identify the space of all analytic extensions of the functions in $\mathcal{R}$ in terms of spaces of analytic functions with some structure.\\

For $0<\alpha \leq 2$ we introduce the open sector $$\Delta_\alpha:=\left\{z\in \mathbb{C}: 0<|\arg(z)|<\pi\alpha/4 \right\}.$$

In \cite{orsini} the author proves that the null reachable space $\mathcal{R}$ is the sum of two Bergman spaces defined on different sectors
$$\mathcal{R}=A^2(\Delta_1)+A^2(1-\Delta_1).$$

In \cite{kellay} the authors improve the last result, they decompose $\mathcal{R}$ as a sum of weighted Bergman spaces
$$\mathcal{R}=A^2(\Delta_1,\omega_{0,\delta})+A^2(1-\Delta_1,\omega_{1,\delta}) \text{ for all }\delta>0,$$
where
$$\omega_{0,\delta}(s)= \delta^{-1} e^{s^2/(2t)}, \,\, s\in \Delta_1, \quad    \omega_{1,\delta}(s)= \delta^{-1}e^{(1-s)^2/(2t)}, \,\, s\in 1-\Delta_1.$$

In \cite{lopez} the author uses the characterization of the image of a certain kind of linear mappings as a RKHS (see Theorem A below or the seminal work \cite{seminal}) to show that $\mathcal{R}$ is a RKHS on $Q=\Delta_1 \cap (1-\Delta_1)$, and computes its reproducing kernel. \\

In this work we study the following 1D degenerate parabolic equation
\begin{eqnarray}\label{sinbeta}
\partial_t u- x^{2(\alpha -1)/\alpha} \partial_{xx} u=0, && x,t>0,  \notag  \\ 
u(0,t)=t^\alpha g(t), && t>0 ,\\
u(x,0)=0, && x>0,\notag\\
\lim_{x\rightarrow \infty}u(x,t)=0, && t>0. \notag
\end{eqnarray}
The next result shows that the solution $u$ of the last system is obtained as a convolution of the data $g$ with a certain positive kernel.
\begin{theorem}\label{resolviendo}
Let $\alpha>0$, $\Rp=\{x\in \mathbb{R}:x>0\}$.   If $u\in C^2(\Rp \times \Rp) \cap C(\overline{\Rp}\times \Rp)$ is a bounded function solving the system (\ref{sinbeta}), then $u(x,t)=t^\alpha\mathcal{L}_t^{\alpha}g (x)$, $x,t >0$, where  
\begin{equation}\label{operaintro}
\mathcal{L}_t^{\alpha}g(x):= \frac{(\alpha/2)^\alpha}{\Gamma(\alpha/2)  t^{\alpha}} \int_0^t \frac{x}{(t-\tau)^{\alpha/2 +1}} \exp\left(- \frac{\alpha^2x^{2/\alpha}}{4(t-\tau)} \right)  g(\tau) \tau^{\alpha}\,d\tau.
\end{equation}
\end{theorem}
\bigskip

Now,  for $t>0$ fixed, we want to describe the space of the analytic extensions $\mathcal{L}_t^{\alpha}g(z)$ with $g$ in a suitable space. Thus, we will prove that $\mathcal{L}_t^{\alpha}$ is defined on a weighted Lebesgue space into a suitable space of analytic function on $\Delta_\alpha$. The case $\alpha=1$ (related to the heat equation) was solved by Aikawa, Hayashi and Saitoh in \cite[Theorem 2.1]{saitoh} and was used as an important step to get the different characterizations of $\mathcal{R}$, see \cite[Lemma 2.5]{hart}, \cite[Proof of Theorem 6]{lopez}, \cite[Proof of Theorem 1.1]{orsini}.\\

For $0<\al \leq 2$ consider the following weighted Bergman space
$$A_{\alpha -1}^2(\Delta_\alpha):=\left\{ f\in hol(\Delta_\alpha): \int_{\Delta_\alpha} |f(z)|^2 (\Re (z^{2/\alpha}))^{\alpha -1}|z|^{(\alpha-2)(\alpha-1)/\alpha}dA(z)<\infty \right \}$$
with the inner product 
$$\langle f,g \rangle_{A_{\alpha -1}^2(\Delta_\alpha)}:=  (\alpha\pi^{1/2})^{\alpha -1}\int_{\Delta_\alpha} f(z) \overline{g(z)} (\Re (z^{2/\alpha}))^{\alpha -1}|z|^{(\alpha-2)(\alpha-1)/\alpha}dA(z),$$
where $dA(z)$ is the Lebesgue measure on $\Delta_\alpha$, and $\Re z$ denotes the real part of $z$.\\

From \cite[Corollary 2.5]{saitoh} we have that 
$$\mathcal{G}_\alpha^t:=  z^{(\alpha -1)(\alpha+2)/(2\alpha)}e^{-\frac{\alpha^2}{4t}z^{2/\alpha}}A^2_{\alpha -1}(\Delta_\alpha)$$
with the norm
\begin{equation}\label{normagalfa}
\|F\|_{\mathcal{G}_\alpha^t}:= \frac{B(\alpha/2,\alpha/2)^{1/2}}{(\al/2)^{(\al -1)/2}\pi^{\frac{\alpha+1}{4}}}  \|z^{(1-\alpha)(\alpha+2)/(2\alpha)}e^{\frac{\alpha^2}{4t}z^{2/\alpha}}F\|_{A^2_{\alpha -1}(\Delta_\alpha)}
\end{equation}
is a RKHS on $\Delta_\al$ with reproducing kernel
$$\frac{(\al/2)^{\al -1}\pi^{\frac{\alpha+1}{2}}}{B(\alpha/2,\alpha/2)}(z\overline{w})^{(\alpha -1)(\alpha+2)/(2\alpha)} e^{-\frac{\alpha^2}{4t}(z^{2/\alpha}+\overline{w}^{2/\alpha} )} K_{\Delta_\alpha,\alpha-1}(z,w),$$
where $K_{\Delta_\alpha,\alpha-1}(z,w)$ is the reproducing kernel of $A^2_{\alpha -1}(\Delta_\alpha)$, see (\ref{meromero}), and $B$ is the Beta function.\\

The main result is the following,
\begin{theorem}\label{main}
For each $t>0$ fixed, the linear mapping
$$\mathcal{L}_t^{\alpha}:L^2_{\mathbb{C}}((0,t),\tau^\alpha d\tau /t^{\alpha})\rightarrow \mathcal{G}_\alpha^t$$ 
is an isometric isomorphism whenever $0<\alpha \leq 2$.
Furthermore, we have the inverse formula
$$(\mathcal{L}_t^{\alpha})^{-1}F(\tau)=\frac{(\alpha/2)^\alpha}{\Gamma(\alpha/2)(t-\tau)^{\alpha/2 +1}}\lim_{N\rightarrow \infty}\int_{E_N}\overline{z}F(z) \exp\left(- \frac{\alpha^2\overline{z}^{2/\alpha}}{4(t-\tau)} \right)d\mu_{\al}^t(z),$$
for all $F\in \mathcal{G}_\alpha^t$ in the topology of $L^2((0,t),\tau^\alpha d\tau/t^{\alpha})$, where $\{E_{N}\}_{N=1}^{\infty}$  is a compact exhaustion of $\Delta_\al$, $\Gamma$ is the Gamma function and
\begin{equation}\label{medidalfa}
d\mu_{\al}^t(z):=\pi^{-1}2^{\al -1}B(\al/2,\al/2) |z|^{4(1-\al)/\al}e^{\frac{\al^2}{2t}\Re(z^{2/\al})}(\Re(z^{2/\al}))^{\al -1}dA(z),
\end{equation}
for $z\in \Delta_\al$.
\end{theorem}
\bigskip

When $ g$ is a continuous function, the next result shows that $u(x,t)=t^\al\mathcal{L}_t^{\alpha}g(x)$ is a classical solution of the 1D degenerate parabolic equation for the half line.
\begin{theorem}\label{nece}
Let $\alpha,T>0$. If $g\in C([0,T])$ then $u(x,t):=t^\alpha\mathcal{L}_t^{\alpha}g(x)\in C^2(\Rp \times (0,T))$ and satisfies 
 \begin{eqnarray}\label{prinsinb}
\partial_t u-  x^{2(\alpha -1)/\alpha} \partial_{xx} u=0, && x>0, \, 0<t<T,  \notag  \\ 
u(0,t)=t^\alpha g(t), && 0<t<T,\\
u(x,0)=0, && x>0,\notag \\
\lim_{x\rightarrow \infty}u(x,t)=0, && 0<t<T. \notag
\end{eqnarray}
\end{theorem}

When $t^\al g$ is a bounded continuous function on $[0, \infty)$ the last result holds for $T= \infty$, see Remark \ref{grandet}.\\

Thus, for $T>0$ fixed we say that $\mathcal{G}_\alpha^T$ is the null reachable space at time $T$ of the degenerate parabolic system (\ref{prinsinb}).\\

This paper is organized as follows. In the next section we consider some results about (weighted) Bergman spaces, and we also include two theorems about RKHS, which are the core of the main result. In Section \ref{degeneratepar} we prove Theorems \ref{resolviendo} and \ref{nece}, and study some properties of the solution to the degenerate parabolic equation in (\ref{prinsinb}). In Section \ref{last} we prove Theorem \ref{main} and give an application.

\section{Preliminaries}
For an open set $\Omega\subset \mathbb{C}$ we denote by $A^2(\Omega)$ the Bergman space on $\Omega$ given by
$$A^2(\Omega):=\left\{f\in hol (\Omega):\int_\Omega |f(z)|^2dA(z) <\infty \right\},$$
and $K_\Omega(z,w)$ stands for the reproducing kernel (the so-called Bergman kernel) of $A^2(\Omega)$.\\

It is well known the conformal invariance of the Bergman kernel: Let $\Omega_1,\Omega_2\subset \mathbb{C}$ be open sets and $\Phi:\Omega_1\rightarrow\Omega_2$ a biholomorphism, then
\begin{equation}\label{kernopeso}
K_{\Omega_1}(z,w)=  \Phi'(z)K_{\Omega_2}(\Phi(z),\Phi(w)) \overline{\Phi'(w)}.
\end{equation}

Now, for $\nu>-1$ we consider the weighted Bergman space $A^2_\nu(\Omega)$ given by
$$A^2_\nu(\Omega):=\left\{ f\in hol(\Omega): \int_\Omega |f(z)|^2K_\Omega(z,z)^{-\nu/2}dA(z) <\infty \right\},$$
and $K_{\Omega,\nu}(z,w)$ denotes its corresponding reproducing kernel (the so-called weighted Bergman kernel).\\

This kind of weighted Bergman kernels are also conformally invariant: Let $\Phi:\Omega_1\rightarrow\Omega_2$ be a biholomorphism such that $\log(\Phi')$ is a well defined holomorphic funtion on $\Omega_1$, then (see \cite[Corollary 6.22]{peloso})
\begin{equation}\label{kerpeso}
K_{\Omega_1,\nu}(z,w)=  \Phi'(z)^{1+\frac{\nu}{2}}K_{\Omega_2,\nu}(\Phi(z),\Phi(w)) \overline{\Phi'(w)^{1+\frac{\nu}{2}}}.
\end{equation}

For instance, the Bergman space on the half plane $\mathbb{C}_+=\{ z \in \mathbb{C}:\Re z >0 \}$ has the Bergman kernel 
$$K_{\mathbb{C}_+}(z,w)=\frac{1}{\pi(z+\overline{w})^2},$$
thus (\ref{kernopeso}) with $\Phi(z)=z^{2/\alpha}$, $0<\alpha\leq 2$, implies that
$$K_{\Delta_\alpha}(z,w)= \frac{4 \,(z\overline{w})^{\frac{2-\alpha}{\alpha}}}{\alpha^2 \pi \left(z^{2/\alpha}+\overline{w}^{2/\alpha}\right)^2}$$
is the Bergman kernel of $A^2(\Delta_\alpha)$.\\
It is well known that (see \cite[Proposition 6.20]{peloso})
$$K_{\mathbb{C}_+\!,\nu}(z,w)=\frac{\nu +1}{\pi^{1+\frac{\nu}{2}}(z+\overline{w})^{\nu+2}},$$
thus (\ref{kerpeso}) implies that 
\begin{equation}\label{meromero}
K_{\Delta_\alpha,\nu}(z,w)=\frac{2^{2+\nu} (\nu+1)\,(z\overline{w})^{\frac{2-\alpha}{\alpha}(1+\frac{\nu}{2})}}{\alpha^{2+\nu} \pi^{1+\frac{\nu}{2}} \left(z^{2/\alpha}+\overline{w}^{2/\alpha}\right)^{\nu+2}}
\end{equation}
is the Bergman kernel of $A^2_\nu(\Delta_\alpha)$.\\

In order to prove the main theorem we introduce a machinery that shows the image of a suitable linear mapping as a RKHS. 
Let $\mathcal{F}(E)$ be the vector space consisting of all complex-valued functions on a set $E$, and let $(\mathcal{H},\langle \cdot, \cdot\rangle_{\mathcal{H}})$ be a Hilbert space. For a mapping $\mathbf{h}:E\rightarrow \mathcal{H}$, consider the induced linear mapping $\mathbf{L}:\mathcal{H}\rightarrow \mathcal{F}(E)$ defined by
$$\mathbf{L}\mathbf{f}(p)=   \langle \mathbf{f}, \mathbf{h}(p) \rangle _{\mathcal{H}}.$$

The vector space $\mathcal{R}(\mathbf{L}):=\{\mathbf{L}\mathbf{f}:\mathbf{f}\in \mathcal{H}\}$ is endowed with the norm
$$\|f\|_{\mathcal{R}(\mathbf{L})}=\inf\{\|\mathbf{f}\|_{\mathcal{H}}:\mathbf{f}\in \mathcal{H}, f=\mathbf{L}(\mathbf{f})\}.$$
A fundamental problem about the linear mapping $\mathbf{L}$ is to characterize the vector space $\mathcal{R}(\mathbf{L})$. The following result summarizes Theorems 2.36, 2.37 in \cite[pages 135--137]{saitohnuevo} and provides an answer to the last question.

\begin{theorem*}
\begin{enumerate}
\item \label{rkhs} $(\mathcal{R}(\mathbf{L}),\|\cdot\|_{\mathcal{R}(\mathbf{L})})$ is a RKHS on $E$ with reproducing kernel
$$\mathbf{K}(p,q)=\langle \mathbf{h}(q),\mathbf{h}(p) \rangle_{\mathcal{H}}, \quad p,q \in E.$$
\item \label{iso}The linear mapping $\mathbf{L}:\mathcal{H}\rightarrow \mathcal{R}(\mathbf{L})$ is an isometric isomorphism iff the set $\{\mathbf{h}(p):p\in E\}$ is complete in $\mathcal{H}$.
\end{enumerate}
\end{theorem*}

We set $H_K(E):=\mathcal{R}(\mathbf{L})$. We assume that $e:H_K(E)\rightarrow L^2(E,du)$ is a continuous embedding and $\mathcal{H}=L^2(I,dm)$, where $(I,dm)$ and $(E,d\mu)$ are $\sigma$-finite measures.\\
Let $h:I \times E \rightarrow \mathbb{C}$ be the function given by
$$h(\tau,p):=\mathbf{h}(p)(\tau), \quad \tau\in I, p\in E.$$
For completeness we reproduce Theorem 2.47 in \cite{saitohnuevo}.
\begin{teorema}
Assume that $\{E_N \}_{N=1}^\infty$ is an increasing sequence of measurable subsets in $E$ such that
$$\bigcup_{N=1}^\infty E_N=E\quad \text{and } \iint_{I\times E_N}|h(\tau,p)|^2dm(\tau)d\mu(p) < \infty \text { for all }N\geq 1.$$ 
Then we have 
$$(e\circ \mathbf{L})^*(f)(\tau)=\lim_{N\rightarrow \infty}\int_{E_N} f(p)h(\tau,p)d\mu(p)$$
for all $f\in L^2(E,d\mu)$ in the topology of $L^2(I,dm)$.
\end{teorema}

\section{On the degenerate parabolic equation}\label{degeneratepar}
In this section we analyze a system very similar to (\ref{prinsinb}). For $\alpha, T>0$ consider the system
\begin{eqnarray}\label{basico}
\partial_t u=  x^{2(\alpha -1)/\alpha} \partial_{xx} u, && x>0, \, 0<t<T,  \notag  \\ 
u(0,t)=g(t), && 0<t<T , \\ \notag
u(x,0)=0, && x>0, \\ \notag
\lim_{x\rightarrow \infty}u(x,t)=0, && 0<t<T. \notag
\end{eqnarray}
In order to get a fundamental solution to the last PDE we follow the ideas when solving the heat equation. For $\alpha>0$ we introduce the generalized complementary (Gaussian) error function
$$E_{\alpha}(\lambda )=\text{erfc}_\alpha(\lambda):=\frac{2}{\Gamma(\alpha/2)}\int_\lambda^\infty \rho^{\alpha -1}e^{-\rho^2}d\rho, \quad \lambda >0.$$
Since there exists a constant $C_\gamma >0$ such that 
\begin{equation}\label{expneg}
e^{-x}\leq C_\gamma x^{-\gamma}
\end{equation}
for all $x,\gamma >0$, we have that $E_\alpha \in C^{\infty}(\mathbb{R}_+)$, and satisfies the following differential equation
\begin{eqnarray}\label{ordinary}
y''(\lambda)+\left(2\lambda -\frac{\alpha -1}{\lambda} \right) y'(\lambda)=0,&& \lambda >0,\\
y(0+)=1, \quad y(\infty)=0. \notag &&
\end{eqnarray}

For $\alpha, x,t >0$ we introduce the function
\begin{eqnarray*}
W_{\alpha}(x,t)&:=& E_{\alpha}\left(  \frac{\al x^{1/\alpha}}{2t^{1/2}} \right) \\
&=& \frac{ \al^\alpha}{2^\al\Gamma(\alpha/2)}\int_0^t \frac{x}{(t-\tau)^{\alpha/2 +1}} \exp\left(- \frac{\al^2 x^{2/\alpha}}{4(t-\tau)} \right)\,d\tau. \notag
\end{eqnarray*}
Since $E_\alpha \in C^{\infty}(\mathbb{R}_+)$ we have that $W_{\alpha}\in C^{\infty}(\mathbb{R}_+\times \mathbb{R}_+)$. By using that $E_\alpha$ satisfies the ODE (\ref{ordinary}), and after some computations, we get that  $W_{\alpha}$ satisfies the PDE in system (\ref{basico}). Moreover,
\begin{equation*}
\lim_{x\rightarrow 0^+}W_{\alpha}(x,t)=1, \lim_{x\rightarrow \infty} W_{\alpha}(x,t)=0 \,\,\forall t>0, \lim_{t\rightarrow 0^+}W_{\alpha}(x,t)=0 \,\, \forall x>0.
\end{equation*}
Notice that $W_{\alpha}$ is the convolution of the constant function 1 with the following positive kernel defined on  $\mathbb{R}_+\times \mathbb{R}_+$,
\begin{equation}
K_{\alpha}(x,t):=\frac{ \al^\alpha}{2^\al\Gamma(\alpha/2)}\frac{x}{t^{\alpha/2 +1}} \exp\left(- \frac{\al^2x^{2/\alpha}}{4t} \right), \quad x,t>0.
\end{equation}

The properties of the function $W_{\alpha}$ suggest to consider the integral transform
\begin{equation}
(\mathcal{T}_t^{\alpha}\,g)(x)= \int_0^t K_{\alpha}(x,t-\tau) g(\tau)\,d\tau, \quad x,t >0,
\end{equation}
where $g$ is a measurable function, in order to get solutions of the system (\ref{basico}).\\

For $\al>0$ we introduce the operator
\begin{equation}
D_{x}^{\alpha}w:=x^{2(\alpha -1)/\alpha} \partial_{xx} w, \quad x>0,
\end{equation}
where $w$ is a function with sufficient regularity.\\

The next result gives some properties of the kernel $K_{\alpha}$ and provides a solution to the system (\ref{basico}) when $g$ is a continuous function.
\begin{proposition}\label{propiedades} Let $\al, T>0$. The following properties hold,
\begin{enumerate}
\item $K_{\alpha}\in C^{\infty}(\mathbb{R}_+\times \mathbb{R}_+)$ and satisfies 
\begin{equation}\label{deriK}
\partial_t K_{\alpha}  =D_{x}^{\alpha} K_{\alpha}.
\end{equation}
In particular,
\begin{equation}\label{deriKgen}
\partial_t^jK_{\alpha}  =(D_{x}^{\alpha})^j K_{\alpha}, \quad j\geq 1.
\end{equation}
\item \label{solucion} If $g\in C([0,T])$ then $u(x,t)=\mathcal{T}_t^{\alpha}\,g(x)$ is a solution of system (\ref{basico}).
\item If $g\in C^m([0,T])$, $m\geq 0$ and $u(x,t)=\mathcal{T}_t^{\alpha}\,g(x)$, then
$$\lim_{x\rightarrow 0+}(D_{x}^{\alpha})^mu(x,t)=\partial_t^mg(t),\quad 0<t<T.$$
\end{enumerate}
\end{proposition}
\begin{proof}  \textit{(1)} A simple computation shows that $\partial_t W_\al = K_\al$, thus $K_\al$ is infinitely differentiable. Since $W_{\alpha}$ satisfies the PDE in system (\ref{basico}) we have
$$\partial_t K_{\alpha}= \partial_t D_{x}^{\alpha}W_\al= D_{x}^{\alpha}\partial_t W_\al=D_{x}^{\alpha} K_\al.$$
We proceed by induction to show that (\ref{deriKgen}) holds. By (\ref{deriK}) the result is valid for $j=1$. Assume that (\ref{deriKgen}) holds for some $j\geq 1$. Since $K_{\alpha}$ is an infinitely differentiable function we have
\begin{eqnarray*}
\partial_t^{j+1}K_{\alpha}  &=& \partial_t^{j} D_{x}^{\alpha} K_{\alpha}= \partial_t^{j}  \left( x^{2(\alpha -1)/\alpha} \partial_{xx}  K_{\alpha} \right)\\
                                                     &=& D_{x}^{\alpha}\partial_t^{j} K_{\alpha}=D_{x}^{\alpha}(D_{x}^{\alpha})^{j} K_{\alpha}=(D_{x}^{\alpha})^{j+1} K_{\alpha}.
\end{eqnarray*}

\textit{(2) }  From (\ref{expneg}) we get
$$\lim_{\tau\rightarrow t^{-}}K_{\alpha}(x,t-\tau)g(\tau)=0,\quad \text{for all }0<t<T.$$
therefore
\begin{eqnarray*}
\partial_t (\mathcal{T}_t^{\alpha}\,g)(x)&=&\int_0^t (\partial_t K_{\alpha})(x,t-\tau) g(\tau)\,d\tau=\int_0^t (D_{x}^{\alpha} K_{\alpha})(x,t-\tau) g(\tau)\,d\tau\\
&=&D_{x}^{\alpha}(\mathcal{T}_t^{\alpha}\,g)(x), \quad \text{for all }x>0, 0<t<T.
\end{eqnarray*}

We make the change of variable 
\begin{equation}\label{cambio}
\rho=\rho(\tau)= \frac{\al x^{1/\alpha}}{2(t-\tau)^{1/2}},
\end{equation}
thus for $x>0, 0<t<T$ we have
\begin{equation}\label{descripcion}
u(x,t)=(\mathcal{T}_t^{\alpha}\,g)(x)=\frac{2}{\Gamma(\alpha /2)}\int_{\frac{\al x^{1/\alpha}}{2t^{1/2}}}^\infty \rho^{\alpha -1}e^{-\rho^2}g\left(t- \frac{\al^2 x^{2/\alpha}}{4\rho^2}\right)d\rho.
\end{equation}
Since $g$ is a continuous function on $[0,T]$, the dominated convergence theorem implies that $u=\mathcal{T}_t^{\alpha}\,g$ satisfies the boundary conditions in (\ref{basico}).\\

\textit{(3)} From (\ref{deriK}) we notice that
$$-\partial_\tau \left(  K_{\alpha}(x,t-\tau) \right)=x^{2(\alpha -1)/\alpha} \partial_{xx}\left( K_{\alpha}(x,t-\tau) \right),$$
therefore
\begin{eqnarray*}
D_{x}^{\alpha}(\mathcal{T}_t^{\alpha}\,g)(x)&=&-\int_0^t g(\tau) \partial_\tau \left(  K_{\alpha}(x,t-\tau) \right) d\tau \\
                               &=& \int_0^t \partial_\tau g(\tau) \,  K_{\alpha}(x,t-\tau) d\tau - g(\tau) K_{\alpha}(x,t-\tau)|_{\tau=0}^{\tau=t}\\
                               &=& (\mathcal{T}_t^{\alpha}\,\partial_t g )(x)+g(0)K_{\alpha}(x,t).
\end{eqnarray*}
By iterating the relation $D_{x}^{\alpha}(\mathcal{T}_t^{\alpha}\,g) =(\mathcal{T}_t^{\alpha}\,\partial_t g )+g(0)K_{\alpha}$ and using (\ref{deriKgen}) we get
\begin{eqnarray*}
(D_{x}^{\alpha})^m(\mathcal{T}_t^{\alpha}\,g) &=&\mathcal{T}_t^{\alpha}\,\partial_t^m g +\sum_{j=0}^{m-1} (\partial_t^{m-1-j} g)(0)(D_x^\alpha)^jK_{\alpha} \\
                                                     &=& \mathcal{T}_t^{\alpha}\,\partial_t^m g +\sum_{j=0}^{m-1} (\partial_t^{m-1-j} g)(0)\partial_t^jK_{\alpha}
\end{eqnarray*}
Using the result in the last item we get
\begin{eqnarray*}
\lim_{x\rightarrow 0+}(D_{x}^{\alpha})^mu(x,t)&=&\lim_{x\rightarrow 0+} (\mathcal{T}_t^{\alpha}\,\partial_t^m g )(x)+\sum_{j=0}^{m-1} (\partial_t^{m-1-j} g)(0)\lim_{x\rightarrow 0+}\partial_t^jK_{\alpha}(x,t) \\
                                        &=& \partial_t^m g(t)+\frac{ (\al/2)^\alpha}{\Gamma(\alpha/2)}\sum_{j=0}^{m-1} (\partial_t^{m-1-j} g)(0)\lim_{x\rightarrow 0+}x\partial_t^j \left[\frac{e^{- \al^2\frac{x^{2/\alpha}}{4t}}}{t^{\alpha /2 +1}} \right]\\
                                        &=& \partial_t^m g(t).
\end{eqnarray*}
\end{proof}

\begin{remark}\label{grandet} By (\ref{descripcion}) it follows that the statement in Proposition \ref{propiedades}-\textit{(\ref{solucion})} holds for $T=\infty$ provided that $g$ is a bounded continuous function on $[0,\infty).$
\end{remark}

\begin{proof}[\textbf{Proof of Theorem \ref{nece}}]
Notice that $t^\al \mathcal{L}^{\alpha}_t g= \mathcal{T}_t^{\alpha}(t^\al g)$ and the result follows by Proposition \ref{propiedades}-\textit{(\ref{solucion})}.
\end{proof}

\begin{proposition} Let $\alpha>0$.   If $u\in C^2(\Rp \times \Rp) \cap C(\overline{\Rp}\times \Rp)$ is a bounded function satisfying 
\begin{eqnarray}\label{sinalfa}
\partial_t u=  x^{2(\alpha -1)/\alpha} \partial_{xx} u, && x,t>0,  \notag  \\ 
u(0,t)=g(t), && t>0 , \\ \notag
u(x,0)=0, && x>0, \\ \notag
\lim_{x\rightarrow \infty}u(x,t)=0, && t>0. \notag
\end{eqnarray}
then $u(x,t)=\mathcal{T}_t^{\alpha}g (x)$, $x,t >0$.
\end{proposition}
\begin{proof} 
We denote by $U (x, s)$ the Laplace transform (denoted by $\mathcal{L}$) of $u$ with respect to the variable $t$, i.e. $U(x,\cdot)=\mathcal{L}(u(x,\cdot)) $. From (\ref{sinalfa}) we get that $U$ satisfies
\begin{eqnarray}\label{ode}
 sU-x^{2(\al-1)/ \alpha}\partial_{xx}U=0&& x>0, \, s>0,  \notag  \\ 
U(0,s)=\mathcal{L} (g), && x>0, \\
\lim_{x\rightarrow \infty}U(x,s)=0, && s>0. \notag 
\end{eqnarray}
We fix $s>0$ and solve the last ODE with respect to the variable $x$.\\

It is well known that the modified Bessel functions of the first kind $I_\nu$ and the second kind $\mathbf{K}_\nu$ are the two linearly independent solutions to the modified Bessel's equation (see \cite[page 374]{abram}):
$$x^2y''+xy'-(x^2+\nu^2)y=0.$$
Moreover (see \cite[page 374]{abram})
$$\lim_{x\rightarrow \infty}\mathbf{K}_\nu (x)=0, \quad \lim_{x\rightarrow \infty}I_\nu(x)=\infty.$$
The Bessel's equation yields
$$s\mathbf{K}_{\alpha /2}''(\alpha s^{1/2}x^{1/\alpha})+\frac{s^{1/2}}{\alpha x^{1/\alpha}}\mathbf{K}_{\alpha /2}'(\alpha s^{1/2}x^{1/\alpha})-\left(s+\frac{1}{4x^{2/\alpha}}\right)\mathbf{K}_{\alpha /2}(\alpha s^{1/2}x^{1/\alpha})=0,$$
which implies the function 
$$V(x,s)=\frac{2\al^{\alpha /2}}{2^{\alpha /2}\Gamma(\alpha/2)} s^{{\alpha/4}}x^{1/2}\mathbf{K}_{\alpha /2}(\alpha s^{1/2}x^{1/\alpha})$$ 
satisfies the ODE and the last condition in (\ref{ode}).\\

Hence there exists a constant $c$ such that $U=cV$. At the end of this proof we will show that
\begin{equation}\label{uno}
\lim_{x\rightarrow 0^+ }V(x,s)=1,
\end{equation}
together with the Dirichlet condition in (\ref{ode}), we get that
$$U=\mathcal{L} (g) V.$$

Since (see \cite[page 917, 8.432-6]{grad})
$$\mathcal{L}\left(\frac{e^{-b/t}}{t^{\nu+1}}\right)=2\left(\frac{s}{b} \right)^{\nu/2}\mathbf{K}_\nu(2\sqrt{sb}), \quad b,s>0,$$
we have that $\mathcal{L}(K_{\alpha}(x,\cdot))=V(x,\cdot)$, then
$$\mathcal{L}(u(x,\cdot))=\mathcal{L} (g)\mathcal{L}(K_{\alpha}(x,\cdot)) = \mathcal{L}(K_{\alpha}(x,\cdot)\ast g)$$
thus
$$u(x,t)=(K_{\alpha}(x,\cdot)\ast g)(t)=\int_0^t K_{\alpha}(x,t-\tau) g(\tau)d\tau=\mathcal{T}_t^{\alpha}g (x). $$
Finally, we use the representation (see \cite[page 917, 8.432-5]{grad})
$$\mathbf{K}_\nu(x)=\frac{2^\nu\Gamma(\nu+\frac{1}{2})}{\sqrt{\pi}x^\nu}\int_0^\infty \frac{\cos(xt)}{(1+t^2)^{\nu+\frac{1}{2}}}dt, \quad \nu\geq -1/2, x>0,$$
and the dominated convergence theorem to obtain (\ref{uno}):
$$\lim_{x\rightarrow 0^+ }V(x,s)=\frac{2\Gamma((\al+1)/2)}{\sqrt{\pi}\Gamma(\al/2)}\int_0^\infty \lim_{x\rightarrow 0^+ } \frac{\cos(\al s^{1/2}x^{1/\al}t)}{(1+t^2)^{\frac{\al}{2}+\frac{1}{2}}}dt=1.$$
\end{proof}

 \begin{proof}[\textbf{Proof of Theorem \ref{resolviendo}}]
 Notice that $t^\al \mathcal{L}^{\alpha}_t g= \mathcal{T}_t^{\alpha}(t^\al g)$ and the result follows by the last proposition.
 \end{proof}

\section{Proof of the main result}\label{last}
Clearly, $K_{\alpha}(\cdot,t)$ is an analytic function on $\Delta_\alpha$ for all $t>0$. Now we pick any $z\in \Delta_\alpha$. By (\ref{expneg}) there exists a constant $C=C(\alpha)>0$ such that
\begin{eqnarray} \label{complex}
| K_{\alpha}(z,t) | &\leq &  \frac{ \al^\alpha}{2^\al \Gamma(\alpha/2)}  \frac{|z|}{t^{\alpha/2+1}}\exp\left(- \frac{\al^2 \Re (z^{2/\alpha})}{4t} \right) \notag \\
                                    &\leq & C|z| (\Re (z^{2/\alpha}))^{-\alpha/2-1} \quad \text{ for all }t>0.
\end{eqnarray}
Therefore $K_{\alpha}(z,\cdot)$ is a bounded continuous function on $(0,\infty).$ So the function $\mathcal{L}_t^{\alpha}g$ is well defined on $\Delta_\al$ for each $g\in L^1_{\mathbb{C}}((0,t),\tau^{\alpha}d\tau)$, $t >0$.\\

The next result shows that the image of $L^1_{\mathbb{C}}((0,t),\tau^{\alpha}d\tau)$ under the mapping $\mathcal{L}_t^{\alpha}$ is a subspace of analytic functions on $\Delta_\alpha$.
\begin{proposition}
Let $t>0$, $0< \al \leq 2$. If $g\in L^1((0,t),\tau^{\alpha}d\tau)$, then $\mathcal{L}_t^{\alpha}g\in hol(\Delta_\alpha).$
\end{proposition}
\begin{proof}
From the estimation (\ref{complex}) and the dominated convergence theorem we have that $\mathcal{L}_t^{\alpha}g$ is a continuous function on $\Delta_\alpha$.
Let $\gamma$ a closed, piecewise differentiable curve in $\Delta_\alpha$, so Fubini's theorem implies that
$$\oint_\gamma \mathcal{L}_t^{\alpha}g(z)dz= \frac{1}{t^{\alpha}} \int_0^t \oint_\gamma  K_{\alpha}(z,t-\tau) \, dz \, g(\tau) \tau^{\alpha}d\tau=0.
$$
From Moreras's theorem we see that $\mathcal{L}_t^{\alpha}g$ is an analytic function on $\Delta_\alpha$.
\end{proof}

\begin{proof}[\textbf{Proof of Theorem \ref{main}}]
Let $\mathcal{H}=L^2_{\mathbb{C}}((0,t),\tau^{\alpha }/t^{\alpha}d\tau)$ with the inner product
$$\langle f,g\rangle_{\mathcal{H}}:=\frac{1}{t^{\alpha}} \int_0^t f(\tau)\overline{g(\tau)}\tau^{\alpha}d\tau.$$
Since $K_{\alpha}(z,\cdot)$ is a bounded continuous function on $(0,\infty)$ for each $z\in \Delta_\al$, the function $\mathbf{h}:\Delta_\alpha \rightarrow \mathcal{H}$ given by
$$\mathbf{h}(z):=\overline{K_{\alpha}(z,t-\cdot)}, \quad z\in \Delta_\alpha$$
is well defined.
Notice that
$$\mathcal{L}_t^{\alpha}g(z)=\langle g, \mathbf{h}(z) \rangle_\mathcal{H}, \quad z\in \Delta_\alpha, g\in \mathcal{H}.$$
Theorem A-(\ref{rkhs}) implies that $\mathcal{R}(\mathcal{L}_t^{\alpha})$ is a RKHS on $\Delta_\alpha$ with reproducing kernel 
\begin{eqnarray*}
\mathcal{K}_{\alpha}(z,w;t)&=&\langle \mathbf{h}(w),\mathbf{h}(z)\rangle_{\mathcal{H}}=\int_0^t K_{\alpha}(z,t-\tau) \overline{K_{\alpha}(w,t-\tau)}\frac{\tau^{\alpha}}{t^{\alpha}}d\tau\\
						&=& \left( \frac{ \al^\alpha}{2^\al\Gamma(\alpha/2)}\right)^2\int_0^t \frac{z\overline{w}}{(t-\tau)^{\alpha +2}} \exp\left(- \al^2\frac{z^{2/\alpha}+\overline{w}^{2/\alpha}}{4(t-\tau)} \right)\frac{\tau^{\alpha}}{t^{\alpha}}  d\tau\\
						&=& \left( \frac{ \al^\alpha}{2^\al\Gamma(\alpha/2)}\right)^2\int_{1/t}^\infty z\overline{w}\left(\rho-t^{-1}\right)^{\alpha} \exp\left(- \al^2(z^{2/\alpha}+\overline{w}^{2/\alpha})\rho/4 \right)  d\rho\\
						&=& \frac{\al^{2\alpha}z\overline{w}}{2^{2\al}(\Gamma(\alpha/2))^2 }e^{-\frac{\alpha^2}{4t}(z^{2/\alpha}+\overline{w}^{2/\alpha} )}\int_0^\infty \eta^{\alpha}  e^{-\al^2(z^{2/\alpha}+\overline{w}^{2/\alpha})\eta/4} d\eta\\
				&=& \frac{4\Gamma(\alpha +1)z\overline{w}}{\alpha^2\left(\Gamma(\alpha/2)\right)^2}\frac{   e^{-\frac{\alpha^2}{4t}(z^{2/\alpha}+\overline{w}^{2/\alpha} )}}{\left(z^{2/\alpha}+\overline{w}^{2/\alpha} \right)^{\alpha+1}} = \frac{4z\overline{w}}{\alpha B(\frac{\alpha}{2},\frac{\alpha}{2})}\frac{   e^{-\frac{\alpha^2}{4t}(z^{2/\alpha}+\overline{w}^{2/\alpha} )}}{\left(z^{2/\alpha}+\overline{w}^{2/\alpha} \right)^{\alpha+1}} \notag    \\
&=& \frac{(\al/2)^{\al -1}\pi^{\frac{\alpha+1}{2}}}{B(\alpha/2,\alpha/2)} (z\overline{w})^{(\alpha -1)(\alpha+2)/(2\alpha)} e^{-\frac{\alpha^2}{4t}(z^{2/\alpha}+\overline{w}^{2/\alpha} )} K_{\Delta_\alpha,\alpha-1}(z,w).
\end{eqnarray*}
It follows that $\mathcal{R}(\mathcal{L}_t^{\alpha})=\mathcal{G}_{ \al}^t$.
\\

Now we claim that $\{\mathbf{h}(z):z \in \Delta_\al\}$ is a complete system in $\mathcal{H}$. Assume that $g\in \mathcal{H}$ satisfies
$$\langle g, \mathbf{h}(z) \rangle_\mathcal{H}=0 \text{ for all }z\in \Delta_\alpha,$$
therefore
$$\int_0^t (t-\tau)^{-(\alpha/2+1)} e^{-(t-\tau)^{-1}z} g(\tau)\tau^{\alpha}d\tau=0 \text{ for all }z\in \mathbb{C}_+.$$
Then we make the change of variable $\rho=\rho(\tau):=(t-\tau)^{-1}-t^{-1} $ to get
\begin{equation}{\label{fourier}}
\int_0^\infty \left(\rho+ t^{-1} \right)^{\frac{\alpha}{2}+1} e^{-x\rho}e^{-i y\rho} g\left(\tau(\rho)\right) \left(\tau(\rho)\right)^{\alpha}d\rho=0
\end{equation}
for all $x>0, y\in \mathbb{R}$. Since $\tau=\tau(\rho)$ is a bounded function and 
$$\int_0^\infty |g\left(\tau(\rho)\right)|^2\left(\rho+t^{-1} \right)^{-2} d\rho<\infty,$$
the factor in  (\ref{fourier}) multiplied by $e^{-iy\rho}$ is in $L^2(\mathbb{R}^+)$. The injectivity of the Fourier transform in $L^2(\mathbb{R})$ implies that $g=0$ a.e. on $(0,t)$.\\

Then Theorem A-(\ref{iso}) implies that $\mathcal{L}_t^{\alpha}:L^2_{\mathbb{C}}((0,t),\tau^{\alpha}/t^{\alpha}d\tau) \rightarrow (\mathcal{G}_\alpha^t,\|\cdot\|_{\mathcal{G}_\alpha^t}) $ is an isometric isomorphism.\\

From (\ref{normagalfa}) and (\ref{medidalfa}) we have that $\mathcal{G}_\alpha^t$ inherits the inner product in $L^2(\Delta_\al, d\mu_\al^t)$. In particular $e: \mathcal{G}_\alpha^t \rightarrow L^2(\Delta_\al, d\mu_\al) $ is a continuous embedding.\\

We just apply Theorem B with $(I,dm)=((0,t),\tau^\al/t^\al d\tau)$ and 
$$\mathbf{L}= \mathcal{L}_t^{\alpha}, \quad h(\tau,z)=\overline{K_{\alpha}(z,t-\tau)}, \,\, 0<\tau<t, z\in \Delta_\al.$$
From (\ref{medidalfa}) and (\ref{complex}) we have that 
$$\iint_{(0,t)\times E} |h(\tau,z)|^2\tau^\al d\tau d\mu^t_\al(z) < \infty $$
for any compact set $E\subset \Delta_\al$, therefore
$$(e\circ \mathcal{L}_t^{\alpha})^*(f)(\tau)=\lim_{N\rightarrow \infty}\int_{E_N}f(z)\overline{K_{\alpha}(z,t-\tau)}d\mu_{\al}^t(z)$$
for all $f\in L^2(\Delta_\al, d\mu_\al^t)$ in the topology of $L^2((0,t),\tau^\al/t^\al d\tau)$, where $\{E_N\}_{N=1}^\infty$ is a compact exhaustion of $\Delta_\al$.\\

Since $\mathcal{L}_t^{\alpha}$ is an isometric isomorphism we have $(\mathcal{L}_t^{\alpha})^*=(\mathcal{L}_t^{\alpha})^{-1}$, so $(e\circ \mathcal{L}_t^{\alpha})^*(f)=(\mathcal{L}_t^{\alpha})^{-1}(f)$ for all $f\in \mathcal{G}_\al^t$, and the result follows.
\end{proof}

When $\al=2$ we have an interesting case,
\begin{corollary} Let $t>0$ fixed. The linear mapping $\mathcal{L}_t:L^2_{\mathbb{C}}((0,t),\tau^{2}/t^{2}d\tau) \rightarrow ze^{-z/t}A^2_1(\mathbb{C}_+)$ given by
\begin{equation*}
\mathcal{L}_tg(z):= \frac{1}{t^2} \int_0^t \frac{z}{(t-\tau)^{2 }} \exp\left(- \frac{z}{t-\tau} \right)  g(\tau) \tau^2\,d\tau
\end{equation*}
is an isometric isomorphism, where $ze^{-z/t}A^2_1(\mathbb{C}_+)$ is endowed with the norm $\|\cdot\|_{\mathcal{G}_2^t}$.
Moreover we have the inverse formula
$$(\mathcal{L}_t)^{-1}F(\tau)=\frac{1}{(t-\tau)^{2}}\lim_{N\rightarrow \infty}\int_{E_N}\overline{z}F(z) \exp\left(- \frac{\overline{z}}{t-\tau} \right)d\mu^t(z),$$
for all $F\in ze^{-z/t}A^2_1(\mathbb{C}_+)$ in the topology of $L^2((0,t),\tau^2 d\tau/t^2)$, where $\{E_{N}\}_{N=1}^{\infty}$  is a compact exhaustion of $\mathbb{C}_+$, and
\begin{equation*}
d\mu^t(z):=2\pi^{-1} |z|^{-2}e^{\frac{2}{t}\Re(z)}\Re(z)dA(z), \quad z\in \mathbb{C}_+.
\end{equation*}
\end{corollary}

As a consequence of Theorem \ref{main} we get an asymptotic behavior of the functions in the image of $\mathcal{L}_t^{\alpha}$. 
\begin{corollary}
Let $g\in L^2((0,t);\tau^\alpha d\tau)$ and $v(\cdot,t)=\mathcal{L}_t^{\alpha}g$. For any $\xi>0$ and $j\geq 0,$ we have
\begin{equation}
\lim_{t\rightarrow 0^+}t^\alpha \left| \partial^j_\xi\left(  v(\xi,t)\xi^{-\frac{(\alpha -1)(\alpha+2)}{2\alpha}}e^{\frac{\alpha^2}{4t}\xi^{2/\alpha}}   \right)\right|^2=0
\end{equation}
\end{corollary}
\begin{proof} By Theorem \ref{main} and the reproducing property of the weighted Bergman kernel $K_{\Delta_\alpha,\alpha-1}$ of $A^2_{\al -1}(\Delta_\al)$ we have
$$v(\xi,t)\frac{e^{\frac{\alpha^2}{4t}\xi^{2/\alpha}}}{\xi^{\frac{(\alpha -1)(\alpha+2)}{2\alpha}}}\!\!=\!\!\iint_{\Delta_\alpha}\!\!v(z,t)\frac{e^{\frac{\alpha^2}{4t}z^{2/\alpha}}}{z^{\frac{(\alpha -1)(\alpha+2)}{2\alpha}}}K_{\Delta_\alpha,\alpha-1}(\xi,z)K_{\Delta_\alpha}(z,z)^{\frac{1-\alpha}{2}} dA(z).$$
Therefore,
\begin{eqnarray*}
\lim_{t\rightarrow 0^+}t^\alpha \left| \partial^j_\xi\left(  v(\xi,t)\frac{e^{\frac{\alpha^2}{4t}\xi^{2/\alpha}}}{\xi^{\frac{(\alpha -1)(\alpha+2)}{2\alpha}}}   \right)\right|^2 \!\!&\leq &\!\! \lim_{t\rightarrow 0^+} t^\alpha \iint_{\Delta_\alpha}\!\! \left| v(z,t) \frac{e^{\frac{\alpha^2}{4t}z^{2/\alpha}}}{z^{\frac{(\alpha -1)(\alpha+2)}{2\alpha}}} \right|^2 K_{\Delta_\alpha}(z,z)^{\frac{1-\alpha}{2}} dA(z)\\
& & \times \iint_{\Delta_\alpha}\!\! |\partial^j_\xi\left(K_{\Delta_\alpha,\alpha-1}(\xi,z)\right) |^2 K_{\Delta_\alpha}(z,z)^{\frac{1-\alpha}{2}} dA(z)\\
&=& \lim_{t\rightarrow 0^+} \int_0^t |g(\tau)|^2\tau^\alpha d\tau  \left[  \partial^j_\xi \partial^j_z  (K_{\Delta_\alpha,\alpha-1}(\xi,z) ) \right]_{z=\xi}.
\end{eqnarray*}
\end{proof}
The following result provides an asymptotic behavior on $\Rp$ of the functions in $A^2_{\alpha -1}(\Delta_{\alpha})$.
\begin{proposition}
Let $j\geq 0$ and $0<\alpha \leq 2$. For any $f\in A^2_{\alpha -1}(\Delta_{\alpha})$ we have
\begin{equation}
\lim_{x\rightarrow\infty}x^{j+(\alpha+1)/2}\partial_x^jf(x)=0, \quad \text{and }\lim_{x\rightarrow0^+}x^{j+(\alpha+1)/2}\partial_x^jf(x)=0.
\end{equation}
\end{proposition}
\begin{proof}
Let $x>0$ fixed. Applying the Cauchy integral formula to $f$ we get
$$\partial_x^jf(x)=\frac{j!}{2\pi i}\oint_{|z-x|=r}\frac{f(z)}{(z-x)^{j+1}}dz$$
where $0<r<x\sin(\pi \alpha/8)$. Thus,
$$x^{j+1}|\partial_x^jf(x)|\leq\frac{j!(j+2)}{2\sqrt{\pi} (\sin(\pi \alpha/8))^{j+1}}\left\{\iint_{\widetilde{D}_x}|f(z)|^2dA(z)\right\}^{1/2}$$
where $\widetilde{D}_x=\{z\in \mathbb{C}:|z-x|<x\sin(\pi \alpha/8)\}$.\\

If $z\in \widetilde{D}_x$, then $z\in \Delta_{\alpha/2}$ and 
$$\left(1-\frac{1}{\sqrt{2}}\right)x\leq  \left(1-\sin\left(\frac{\pi\alpha}{8}\right)\right)x< |z|<\left(1+\sin\left(\frac{\pi\alpha}{8}\right)\right)x\leq \left(1+\frac{1}{\sqrt{2}}\right)x,$$
therefore
$$\frac{1}{\sqrt{2}}\left(1-\frac{1}{\sqrt{2}}\right)^{2/\alpha}x^{2/\alpha}<\Re(z^{2/\alpha})=|z|^{2/\alpha}\cos(2\arg(z)/\alpha)\leq \left(1+\frac{1}{\sqrt{2}}\right)^{2/\alpha}x^{2/\alpha}$$
for all $z\in \widetilde{D}_x$.\\
The last inequalities imply that
$$\frac{(1-1/\sqrt{2})^{2(2-\alpha)/\alpha}}{(1+1/\sqrt{2})^{4/\alpha}} x^{ -2}   <\alpha^2\pi K_{\Delta_\alpha}(z,z)< 2\frac{(1+1/\sqrt{2})^{2(2-\alpha)/\alpha}}{(1-1/\sqrt{2})^{4/\alpha}}  x^{ -2}$$
for all $z\in \widetilde{D}_x$, it follows that 
$$K_{\Delta_\alpha}(z,z)^{(1-\alpha)/2} \approx x^{\alpha -1} \quad \text{on } \widetilde{D}_x.$$
Hence there exists a constant $C>0$ such that
$$x^{j+(\alpha+1)/2}|\partial_x^jf(x)|\leq C \left\{\iint_{\widetilde{D}_x}|f(z)|^2K_{\Delta_\alpha}(z,z)^{(1-\alpha)/2}dA(z)\right\}^{1/2}$$
for all $x>0$, and the result follows by the dominated convergence theorem.
\end{proof}

\textbf{Conclusion} Recently, Cannarsa et. al. proved a suitable global Carleman estimate  to get the null controllability for a degenerate parabolic equation on a finite interval, see \cite{cannarsa}. As in the heat equation case, now the problem is to characterize the null reachable space of the degenerate parabolic equation studied by them, and this work is a first step to achieve that goal.


\begin{thebibliography}{99}

\bibitem{abram}  Abramowitz, M. and Stegun, I. A. (Eds.). Handbook of Mathematical Functions with Formulas, Graphs, and Mathematical Tables, 9th printing. New York: Dover, 1972.

\bibitem{saitoh} Aikawa, Hiroaki; Hayashi, Nakao; Saitoh, Saburou. The Bergman space on a sector and the heat equation. Complex Variables Theory Appl. 15 (1990), no. 1, 27--36.

\bibitem{cannarsa} P. Cannarsa, P. Martinez, and J. Vancostenoble; Carleman estimates for a class of degenerate parabolic operators, SIAM J. Control Optim. 47 (2008), no. 1, 1–19

\bibitem{grad} Gradshteyn, I. S.; Ryzhik, I. M.; Table of integrals, series, and products. Translated from the Russian. Seventh edition. Elsevier/Academic Press, Amsterdam, 2007. xlviii+1171 pp.

\bibitem{hart} Hartmann Andreas, Kellay Karim, Tucsnak Marius. From the reachable space of the heat equation to Hilbert spaces of holomorphic functions, hal-01569695, 2017.

\bibitem{kellay} Kellay, Karim; Normand, Thomas; Tucsnak, Marius. Sharp reachability results for the heat equation in one space dimension, hal-02302165.

\bibitem{lopez}L\'opez-Garc\'ia, Marcos. The reachable space of the heat equation for a finite rod as a Reproducing Kernel Hilbert Space. arXiv:1910.03765

\bibitem{orsini}  Orsini, Marcu-Antone. Reachable states and holomorphic function spaces for the 1-D heat equation. arXiv:1909.01644.

\bibitem{peloso}Peloso, Marco M. Classical spaces of holomorphic functions. http://www.mat.unimi.it/ users/peloso/Matematica/har-ber2a.pdf

\bibitem{seminal} Saitoh, Saburou; Hilbert spaces induced by Hilbert space valued functions. Proc. Amer. Math. Soc. 89 (1983), no. 1, 74–78.

\bibitem{saitohnuevo} Saitoh, Saburou; Sawano, Yoshihiro. Theory of reproducing kernels and applications. Developments in Mathematics, 44. Springer, Singapore, 2016. xviii+452 pp.


\end{thebibliography}
\end{document}